# $C-$ortocentros y sistemas $C-$ortocéntricos en planos de Minkowski

Tobías de Jesús Rosas Soto
(tjrosas@gmail.com; trosas@demat-fecluz.org)

Departamento de Matemática
Facultad Experimental de Ciencias
Universidad del Zulia
Venezuela

## Resumen

Usando la noción de $C-$ortocentro se extienden, a planos de Minkowski en general, nociones de la geometría clásica relacionadas con un triángulo, como por ejemplo: puntos de Euler, triángulo de Euler, puntos de Poncelet. Se muestran propiedades de estas nociones y sus relaciones con la circunferencia de Feuerbach. Se estudian sistemas $C-$ortocéntricos formados por puntos presentes en dichas nociones y se establecen relaciones con la ortogonalidad isósceles y cordal. Además, se prueba que la imagen homotética de un sistema $C-$ortocéntrico es un sistema $C-$ortocéntrico.

*Palabras claves:* $C-$ortocentro, sistemas $C-$ortocéntricos, planos de Minkowski, puntos de Euler, puntos de Poncelet, Circunferencia de Feuerbach.

## Abstract

Using the notion of $C-$orthocenter, notions of the classic euclidean geometry related with a triangle, as for example: Euler points; Euler's triangle; and Poncelet's points, are extended to Minkowski planes in general. Properties of these notions and their relations with the Feuerbach's circle, are shown. $C-$orthocentric systems formed by points in the above notions are studied and relations with the isosceles and chordal orthogonality, are established. In addition, there is proved that the homothetic image of a $C-$orthocentric system is a $C-$orthocentric system.

*Key words:* $C-$orthocenter, $C-$orthocentric systems, Minkowski planes, Euler points, Poncelet's points Feuerbach circle.

**1. Introducción.**

En 1960, Asplun y Grünbaum introdujeron las nociones de $S_X-$ortocentro y sistemas $S_X-$ortocéntricos en planos de Minkowski estrictamente convexos y suaves, y mostraron algunas propiedades de las mismas (ver [5]). En 2007, la matemática M. Spirova y el matemático H. Martini establecieron que la condición de suavidad podía ser omitida, y mostraron una serie de teoremas de $C-$ortocentricidad relacionados con la circunferencia de Feuerbach y otras circunferencias relacionadas con los triángulos en planos de Minkowski estrictamente convexos (ver [8]).

En el presente año (2014), los matemáticos T. Rosas y W. Pacheco mostraron en [15] que la condición de convexidad estricta también puede ser omitida para definir la noción de $C-$ortocentro en planos de Minkowski en general y por tanto, la de sistemas $C-$ortocéntricos. Esto les permitió generalizar las nociones de recta de Euler y circunferencia de Feuerbach para planos de Minkowski en general.

En este trabajo se presentan extensiones a planos de Minkowski en general de nociones existentes en la geometría euclidiana clásica relacionadas con un triángulo, que se definen en función del ortocentro del mismo, tales como: puntos de Euler, triángulo de Euler, triángulo medial, puntos de Poncelet y circunferencia de Feuerbach. Se muestra que sustituyendo la noción de ortocentro por la de $C-$ortocentro, tales nociones se pueden definir en planos de Minkowski en general y se prueba la validez de algunas propiedades geométricas que se verifican en la geometría euclidiana clásica. Además, se estudian algunas propiedades geométricas de los sistemas $C-$ortocéntricos relacionadas con las nociones mencionadas y las ortogonalidades isósceles y cordal.

A pesar que trabajaremos con planos de Minkowski en general y por tanto, las circunferencias asociadas a la norma del plano pueden ser polígonos, con un número par de lados, se utilizará la circunferencia asociada a la norma euclidia en la mayoría de las figuras geométricas que se presentarán para ilustrar, de mejor manera, las ideas de las demostraciones. Para estudiar la geometría en planos de Minkowski y las propiedades básicas de las ortogonalidades isósceles y cordal véanse las referencias [6, 8, 10, 11, 15, 16] y la monografía [1].

Denotemos por $(\mathbf{R}^2, \|\circ\|) = M$ a un plano de Minkowski cualquiera con origen $O$, circunferencia unitaria $C$ y norma $\|\circ\|$. Para cualquier punto $x \in M$ y $\lambda \in \mathbf{R}^+$, llamemos al conjunto $C(x, \lambda) = x + \lambda C$ la *circunferencia de centro en* $x$ y *radio* $\lambda$. Para $x \neq y$, denotemos por $\langle x, y \rangle$ a la *línea* que pasa por $x$ e $y$, por $[x, y]$ el *segmento* entre $x$ y $y$, y por $[x, y\rangle$ el *rayo* que inicia en $x$ y pasa a través de $y$ (ver [15, 16]).

Para puntos $x, y \in M$, se dice que $x$ es *ortogonal isósceles* a $y$ si $\|x + y\| = \|x - y\|$, y lo denotaremos por $x \perp_I y$. De igual forma, $x$ se dice *ortogonal Birkhoff* a $y$ si $\|x + ty\| \geq \|x\|$ para todo $t \in \mathbf{R}$, y lo denotaremos por $x \perp_B y$ (ver [10, 11]). Por otra parte, sean $L_1 = [p_1, q_1]$ y $L_2 = [p_2, q_2]$ cuerdas de una circunferencia $G$ en $M$, entonces $L_1$ y $L_2$ se dicen *ortogonales cordales* si la línea que pasa a través de $q_2$ y $p'_2$ (el opuesto de $p_2$ en $G$) es paralela a la línea $\langle p_1, q_1 \rangle$. En caso de que $p'_2 = q_2$, se dirá que $L_1$ es *ortogonal cordal* a $L_2$ si existe una línea soporte de $G$ que pase por $q_2$ y sea paralela a la línea $\langle p_1, q_1 \rangle$. Si $q'_2$ es el opuesto de $q_2$ en $G$, entonces $p'_2 q_2 p_2 q'_2$ es un paralelogramo (ver [8]). Consideraremos la ortogonalidad cordal de forma natural solo con respecto a la circunferencia unitaria $C$, usando el símbolo $\perp_C$ para esto, es decir, escribir $[p_1, q_1] \perp_C [p_2, q_2]$ automáticamente presupone que $[p_1, q_1]$ y $[p_2, q_2]$ son cuerdas de $C$. Algunas propiedades necesarias de estas ortogonalidades son las siguientes (ver [8, 10, 11]): para $x, y \in M$ tales que $x \perp_I y$, entonces $y \perp_I x$, es decir, la ortogonalidad isósceles es simétrica; y para toda cuerda $[x, y] \in C$ existe una cuerda $[z, w] \in C$ tal que $[x, y] \perp_C [z, w]$.

Para $p \in M$, denotemos por $S_p$ a la *simetría* con respecto al punto $p$, dada por la expresión $S_p(w) = 2p - w$ para $w \in M$. $H_{p,k}$ denotará la *homotecia* con centro $p$ y razón $k$ ($k \in \mathbf{R}$), y está dada por $H_{p,k}(w) = (1-k)p + kw$ para $w \in M$. Recordemos que las simetrías son *isometrías* en planos de Minkowski, es decir, $\|S_p(w) - S_p(v)\| = \|w - v\|$ para todo $w, v \in M$ (ver [15, 16]). Para $x, z \in M$, sean $e$ un punto en el segmento $[x, z]$, $w$ uno sobre la prolongación del $[x, z]$ y $z_0 = S_w(z)$. Diremos que los puntos $e$ y $w$ son *conjugados armónicos* de $x$ y $z$, si $H_{e,-k}(z) = x$ y $H_{w,-k}(z_0) = x$, con $z \in \mathbf{Z}$ (ver [4, 16]).

Para puntos $x_1, x_2, x_3 \in M$ denotemos por $\Delta x_1 x_2 x_3$ al *triángulo* de vértices $x_1, x_2, x_3$. Si $p_4 \in M$, diremos que el $\Delta p_1, p_2, p_3$ es el $p_4$–antitriángulo del $\Delta x_1 x_2 x_3$, si $p_i = S_{m_i}(p_4)$ para $i = 1, 2, 3$, con $m_i$ los puntos medios de los lados del $\Delta x_1 x_2 x_3$. Al $\Delta m_1 m_2 m_3$ lo llamaremos *triángulo medial* del $\Delta x_1 x_2 x_3$. Diremos que $p_4$ es un *circuncentro* del $\$\backslash\Delta x_1 x_2 x_3$, si $\|p_4 - x_1\| = \|p_4 - x_2\| = \|p_4 - x_3\|$. Si tal $p_4$ existe, es el centro de una circunferencia que pasa por los vértices del $\Delta x_1 x_2 x_3$, que llamaremos *circunferencia circunscrita* y a su radio el *circunradio* (ver [15, 16]). Denotemos por $C(\Delta x_1 x_2 x_3)$ el conjunto de los circuncentros del $\Delta x_1 x_2 x_3$. Si $p_4 \in C(\Delta x_1 x_2 x_3)$, diremos que $x_4$ es el $C$–ortocentro del $\Delta x_1 x_2 x_3$ asociado a $p_4$ si $S_q(p_4) = x_4$, donde $q$ es el punto de simetría del $\Delta x_1 x_2 x_3$ y su $p_4$–antitriángulo (ver [15, 16]).

El conjunto $\{x_1, x_2, x_3, x_4\}$ se dirá sistema $C$–ortocéntrico, si existe $p_4 \in C(\Delta x_1 x_2 x_3)$ tal que $S_q(p_4) = x_4$, donde $q$ es el punto de simetría del $\Delta x_1 x_2 x_3$ y su $p_4$–antitriángulo. Denotemos por $H(\Delta x_1 x_2 x_3)$ el conjunto de $C$–ortocentros del $\Delta x_1 x_2 x_3$. Llamaremos circunferencia de Feuerbach del $\Delta x_1 x_2 x_3$ a la circunferencia que pasa por los puntos medios del $\Delta x_1 x_2 x_3$ y los puntos medios de los segmentos formados por el $C$–ortocentro del $\Delta x_1 x_2 x_3$ y los vértices del mismo. (ver [9, 15, 16])

2. **Preliminares.**

Los siguientes resultados son necesarios para la investigación.

**Lema 2.1:** (ver [8]) Para puntos diferentes $x, y, z \in C$, con $z \neq -x$ y $z \neq -y$, las relaciones $[x,z] \perp_C [y,z]$ y $[y,z] \perp_C [x,z]$ se cumplen si y solo si $x$ e $y$ son puntos opuestos en $C$.

**Teorema 2.1:** (ver [15, 16]) Sea $M$ un plano de Minkowski. Sean $x_1, x_2, x_3$ y $p_4$ puntos en $M$. Sean $m_1, m_2$ y $m_3$ los puntos medios de los segmentos $[x_2, x_3]$, $[x_1, x_3]$ y $[x_1, x_2]$, respectivamente. Definamos los puntos $p_i = S_{m_i}(p_4)$, para $i = 1, 2, 3$, entonces se cumple lo siguiente:

1. Los segmentos $[x_i, p_i]$ tienen el mismo punto medio $q$, para $i = 1, 2, 3$. Además, $2(q - m_i) = x_i - p_4$ para $i = 1, 2, 3$, es decir, $q = \dfrac{x_1 + x_2 + x_3 - p_4}{2}$.
2. Si $x_4 = S_q(p_4)$, entonces $x_i - x_j = p_j - p_i$ para $\{i, j\} \subset \{1, 2, 3, 4\}$.
3. $x_i - p_j = p_k - x_l$, donde $\{i, j, k, l\} = \{1, 2, 3, 4\}$.
4. Si $g = \dfrac{x_1 + x_2 + x_3}{3}$, entonces $H_{g, -2}(p_4) = x_4$.

**Corolario 2.1:** (ver [15, 16]) Con la hipótesis del teorema previo, se cumple lo siguiente:

1. Si $x_4 = S_q(p_4)$, entonces $x_4 = x_1 + x_2 + x_3 - 2p_4$.
2. Si $g = \dfrac{x_1 + x_2 + x_3}{3}$ y $g_1 = \dfrac{p_1 + p_2 + p_3}{3}$, entonces $g_1 = S_q(g)$ y $S_g(p_4) = g_1$.

1. **Resultados y pruebas.**

En esta sección se presentan los resultados principales del trabajo.

**Teorema 3.1:** Sean $p_4 \in C(\Delta x_1 x_2 x_3)$, $\Delta p_1, p_2, p_3$ el $p_4$–antitriángulo del $\Delta x_1 x_2 x_3$, $q$ el punto de simetría de dichos triángulos y $x_4 = S_q(p_4)$, entonces se cumple:

1. Los puntos $p_i$ y $x_i$ son circuncentros de los triángulos $\Delta x_j x_k x_l$ y $\Delta p_j p_k p_l$, para $\{i, j, k, l\} = \{1, 2, 3, 4\}$, respectivamente. Los circunradios de dichos triángulos son iguales entre sí.
2. Los puntos medios de los lados del $\Delta x_j x_k x_l$ y su $p_i$–antitriángulo, están en la circunferencia de centro $q$ y radio $\|q - m_l\|$, para $\{i, j, k, l\} = \{1, 2, 3, 4\}$.
3. $C(q, \|q - m_l\|) = H_{x_4, \frac{1}{2}}(C(p_4, \|p_4 - x_l\|))$. En particular, $q = H_{x_4, \frac{1}{2}}(p_4)$.

*Demostración:* 1. Dado que $p_4 \in C(\Delta x_1 x_2 x_3)$ (ver *Figura* 1) se tiene que $\|p_4 - x_i\| = \alpha$ para $i = 1,2,3$, con $\alpha \in \mathbb{R}$. Por los ítem 2 y 3 del Teorema 2.1, $\|x_4 - p_i\| = \alpha$ y $\|p_j - x_k\| = \alpha$ para $\{i,j,k\} = \{1,2,3\}$. Por tanto, $\|p_i - x_j\| = \alpha$ para $\{i,j\} \subset \{1,2,3,4\}$. Así, $p_i \in C(\Delta x_j x_k x_l)$ y $x_i \in C(\Delta p_j p_k p_l)$ para $\{i,j,k,l\} = \{1,2,3,4\}$.

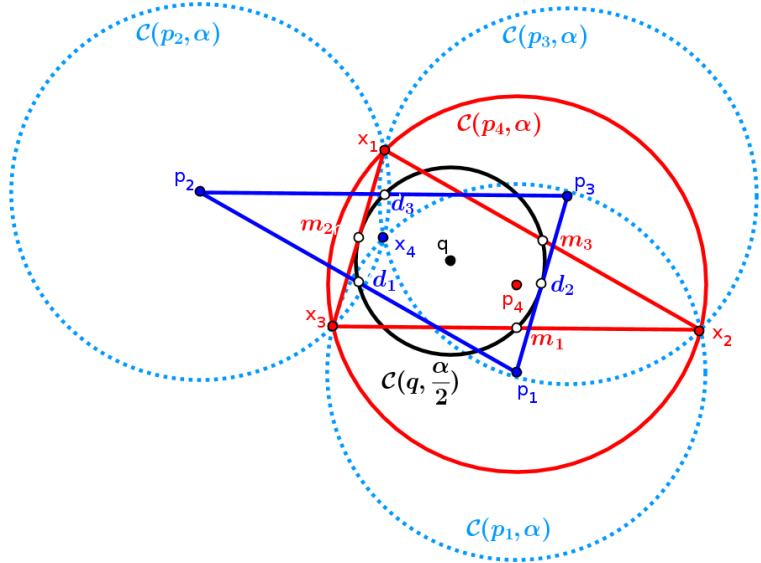

Figura 1: *Circunferencia de los seis puntos* $C(q, \tfrac{\alpha}{2})$

2. Sea $m_{ij} = \dfrac{x_i + x_j}{2}$ y para $\{i,j\} \subset \{1,2,3,4\}$. Veamos que $C(q, \|q - m_l\|)$ pasa por los puntos medios del $\Delta x_j x_k x_l$ (ver *Figura 1*). Como $\|x_i - p_4\| = \alpha$ para $i = 1,2,3$, entonces por el ítem 1 del Teorema 2.1 tenemos que

$$\|q - m_{ij}\| = \left\| \frac{x_1 + x_2 + x_3 - p_4}{2} - \left(\frac{x_i + x_j}{2}\right) \right\| = \left\| \frac{x_k - p_4}{2} \right\| = \frac{\alpha}{2}$$

para $\{i,j,k\} = \{1,2,3\}$. Supongamos que $i = 4$, entonces por el ítem 1 del Teorema 2.1 obtenemos que

$$\|q - m_{j4}\| = \left\| \frac{x_1 + x_2 + x_3 - p_4}{2} - \left(\frac{x_j + x_4}{2}\right) \right\| = \left\| \frac{x_j + x_k}{2} - q \right\| = \|m_{jk} - q\| = \frac{\alpha}{2}$$

para $\{j,k\} \subset \{1,2,3\}$. De igual forma se razona si $j = 4$.

Veamos que $C(q, \|q - m_l\|)$ pasa por los puntos medios del $p_i$-antitriángulo del $\Delta\, x_j x_k x_l$.

Sean $d_{ij} = \dfrac{p_i + p_j}{2}$, para $\{i, j\} \subset \{1,2,3,4\}$, entonces

$$d_{ij} = \frac{p_i + p_j}{2} = \frac{4q - (x_i + x_j)}{2} = 2q - m_{ij}$$

y por tanto $\|q - d_{ij}\| = \|q - m_{ij}\|$ para $\{i, j\} \subset \{1,2,3,4\}$. Así, obtenemos lo deseado.

3. Tomemos $w \in C(p_4, \|p_4 - x_l\|)$, por tanto $\|p_4 - w\| = \|p_4 - x_l\|$. Veamos que $H_{x_4, \frac{1}{2}}(w) \in C(q, \|q - m_l\|)$. Luego, por hipótesis obtenemos que:

$$\left\| q - H_{x_4, \frac{1}{2}}(w) \right\| = \left\| q - \frac{x_4 + w}{2} \right\| = \frac{1}{2} \|p_4 - w\| = \frac{1}{2} \|p_4 - x_l\|$$

Como $x_l - p_4 = 2(q - m_l)$ por el ítem 1 Teorema 2.1, entonces $\left\| q - H_{x_4, \frac{1}{2}}(w) \right\| = \|q - m_l\|$ obteniendo lo deseado. □

Un resultado que se obtiene de forma directa del Teorema 3.1 es el siguiente:

**Corolario 3.1:** Con las hipótesis del teorema anterior se cumple lo siguiente:

1. La circunferencia de los seis puntos (circunferencia de Feuerbach) de los triángulos $\Delta\, p_j p_k p_l$ y $\Delta\, x_j x_k x_l$ coinciden, y su centro es $q$.
2. $q = \dfrac{x_i + x_j + x_k - p_l}{2}$ para $\{i, j, k, l\} = \{1,2,3,4\}$.
3. $x_i = x_j + x_k + x_l - 2 p_i$ para $\{i, j, k, l\} = \{1,2,3,4\}$.

En la geometría euclidiana clásica las nociones de segmentos de Euler y segmentos de Poncelet están intrínsecamente relacionadas con un triángulo, su ortocentro y la circunferencia circunscrita del mismo. Estas se definen como los segmentos formados por un vértice del triángulo y el ortocentro, y el ortocentro y un punto de la circunferencia circunscrita, respectivamente. A los puntos medios de dichos segmentos se les denominan puntos de Euler y puntos de Poncelet, respectivamente.

Como se puede ver en [15, 16], la noción de $C$ – ortocentro en planos de Minkowski coincide con la noción de ortocentro en los planos euclídeos. Por tanto, es natural pensar en definir las estructuras antes mencionadas de la siguiente forma:

**Definición 3.1:** Sean el $\Delta x_1 x_2 x_3$ en un plano de Minkowski $M$ y $x_4 \in H(\Delta x_1 x_2 x_3)$. Se llamará *segmento de Euler*, al segmento $[x_i, x_4]$ para $i = 1,2,3$. Al punto medio de dichos segmentos se les llamará *puntos de Euler*.

De manera similar definiremos las nociones de segmentos y puntos de Poncelet por:

**Definición 3.2:** Sean el $\Delta x_1 x_2 x_3$ en un plano de Minkowski $M$, $x_4 \in H(\Delta x_1 x_2 x_3)$ y $C(p_4, \lambda)$ la circunferencia circunscrita del triángulo, asociada a $x_4$. Se llamará *segmento de Poncelet*, al segmento $[w, x_4]$, con $w \in C(p_4, \lambda)$. Al punto medio de dichos segmentos se les llamará *puntos de Poncelet*.

Tomando en cuenta estas definiciones, el ítem 2 del Corolario 2.1 y el ítem 2 del Teorema 3.1, podemos ver que es valido el siguiente resultado:

**Corolario 3.2:** Dado el $\Delta x_1 x_2 x_3$ y $p_4 \in C(\Delta x_1 x_2 x_3)$, entonces
1. Cada punto de Euler del $\Delta x_1 x_2 x_3$ es un punto de Poncelet de este.
2. Cada punto de Euler del $p_4$ – antitriángulo del $\Delta x_1 x_2 x_3$, es un punto de Poncelet del $p_4$ – antitriángulo.
3. Los puntos de Euler del $\Delta x_1 x_2 x_3$ y su $p_4$ – antitriángulo, están en la circunferencia de los seis puntos del $\Delta x_1 x_2 x_3$.

El siguiente lema muestra ciertas propiedades del triángulo de Euler y el medial de un triángulo dado.

**Lema 3.1:** Sea el $\Delta x_1 x_2 x_3$, $p_4 \in C(\Delta x_1 x_2 x_3)$ y $x_4 \in H(\Delta x_1 x_2 x_3)$ tal que $x_4 = S_q(p_4)$, donde $q$ es punto de simetría entre el triángulo dado y su $p_4$ – antitriángulo, entonces se cumple:

1. El triángulo de Euler y el medial del $\Delta x_1 x_2 x_3$ son congruentes.
2. $q$ es circuncentro del triángulo de Euler y del triángulo medial del $\Delta x_1 x_2 x_3$. Además, el circunradio de dichos triángulos es la mitad del circunradio del $\Delta x_1 x_2 x_3$.

3. Los puntos $p_4$ y $x_4$ son $C$ – ortocentro, relacionados con $q$, del triángulo medial y el triángulo de Euler del $\Delta x_1x_2x_3$, respectivamente.

4. El punto $S_{m_i}(x_4)$ está en la circunferencia circunscrita, para todo $i = 1,2,3$, y es diametralmente opuesto al vértice opuesto al lado, es decir, $S_{m_i}(x_4) = S_{p_4}(x_i)$ para $i = 1,2,3$.

*Demostración:* 1. Sean $d_i = \dfrac{S_{m_j}(p_4) + S_{m_k}(p_4)}{2}$ para $\{i,j,k,l\} = \{1,2,3,4\}$. Se tiene que $d_i - d_j = m_j - m_i$, de manera que los triángulos $\Delta x_1x_2x_3$ y $\Delta p_jp_kp_l$ son congruentes (ver *Figura 2*).

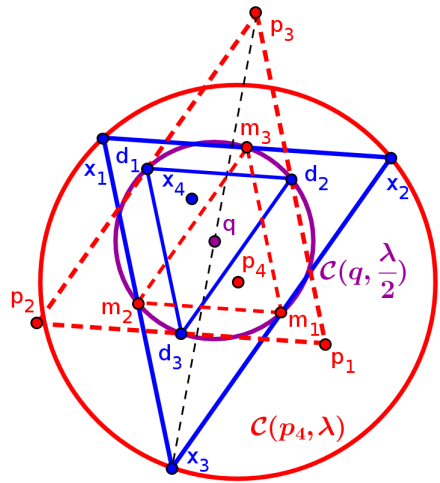

Figura 2: *Demostración de los ítem 1 y 2, Lema 3.1.*

2. Por el ítem 2 del Corolario 2.1, el triángulo medial del $p_4$ – antitriángulo del $\Delta x_1x_2x_3$, es el triángulo de Euler del $\Delta x_1x_2x_3$. Por el ítem 2 del Teorema 3.1, los vértices de dichos triángulos están en la circunferencia $C\left(q, \frac{\|p_4 - x_1\|}{2}\right)$. Por tanto, esta es la circunscrita del triángulo de Euler y del medial asociados al $\Delta x_1x_2x_3$, y su radio es la mitad del circunradio del $\Delta x_1x_2x_3$ (ver *Figura 2*).

3. Sea $\|p_4 - x_i\| = \alpha$ para $i = 1,2,3$, y $w_i = \dfrac{m_j + m_k}{2}$ para $\{i,j,k\} = \{1,2,3\}$. Por el ítem 2 del Teorema 3.1, el punto $q$ es el circuncentro del $\Delta m_1m_2m_3$. Además, $q'$ es el punto de simetría de los triángulos $\Delta m_1m_2m_3$ y $\Delta q_1q_2q_3$, con $q_i = S_{w_i}(q)$ para $i = 1,2,3$. Por tanto $S_{q'}(q) = p_4$, es decir, $p_4$ es $C$ – ortocentro del triángulo medial $\Delta m_1m_2m_3$. Luego, como $d_i = S_q(m_i)$ para $i = 1,2,3$, y $x_4 = S_q(p_4)$, entonces $x_4 \in H(\Delta d_1d_2d_3)$ (ver *Figura 3*).

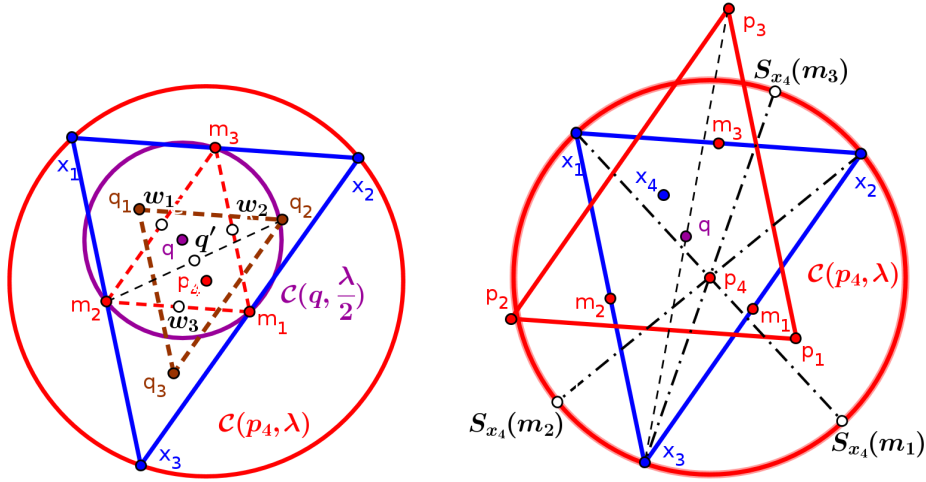

Figura 3: *Demostración de los ítem 3 y 4, Lema 3.1.*

4. Como $p_4$ es circuncentro del $\Delta x_1 x_2 x_3$, entonces $\|p_4 - x_i\| = \lambda$ para $i = 1,2,3$. Luego, $\|p_4 - S_{m_i}(x_4)\| = \|p_4 - 2m_i + x_4\|$ y, por el ítem 1 del Corolario 2.1, se tiene que $\|p_4 - S_{m_i}(x_4)\| = \|x_i - p_4\| = \lambda$, por tanto $S_{m_i}(x_4) \in C(p_4, \lambda)$ para $i = 1,2,3$ y además, $S_{p_4}(x_i) = 2p_4 - x_i = x_j + x_k - x_4 = S_{m_i}(x_4)$ para $\{i, j, k\} = \{1,2,3\}$ (ver *Figura 3*). □

Una de las cosas más relevantes que muestra el siguiente lema es la conjugación armónica existente entre el centro de la circunferencia de Feuerbach y el circuncentro de un triángulo dado, con el $C$–ortocentro y el baricentro del mismo.

**Lema 3.2:** Con las mismas hipótesis del lema precedente, sea $C(p_4, \lambda)$ la circunferencia circunscrita del $\Delta x_1 x_2 x_3$ y $m_i$ el punto medio del segmento $[x_j, x_k]$ para $\{i, j, k\} = \{1,2,3\}$, entonces se cumple:

1. El simétrico del punto medio de un lado del triángulo, con respecto al centro de su circunferencia de los seis puntos, es un punto de Euler, es decir, $S_q(m_i)$ es un punto de Euler para $i = 1,2,3$.
2. El simétrico de los vértices del triángulo, respecto del centro de la circunferencia de los seis puntos, es el simétrico de $p_4$ con respecto al punto medio del lado opuesto, es decir, $S_q(x_i) = S_{m_i}(p_4)$ para $i = 1,2,3$.
3. El centro de la circunferencia de los seis puntos, $q$, y $p_4$ son conjugados armónicos respecto del segmento que une al $C$–ortocentro y el baricentro del $\Delta x_1 x_2 x_3$.

4. Si $\left[S_{p_4}(x_i), x_i\right]$ es un circundiámetro de $C(p_4, \lambda)$ para $i = 1,2,3$, entonces los puntos $S_{p_4}(x_i)$, $m_i$ y $x_4$ son puntos colineales, es decir, $x_4 \in \langle S_{p_4}(x_i), m_i \rangle$ para $i = 1,2,3$.

*Demostración:* 1. Dado que $x_4$ es $C$–ortocentro del $\triangle x_1 x_2 x_3$ asociado a $p_4$, entonces por el ítem 1 del Corolario 2.1, $x_4 = x_1 + x_2 + x_3 - 2p_4$. Por el ítem 1 del Teorema 2.1, se tiene que $S_q(m_i) = x_1 + x_2 + x_3 - p_4 - m_i$ para $i = 1,2,3$, es decir,

$$S_q(m_i) = \frac{x_j + x_k}{2} + x_i - p_4 = \frac{x_j + x_k + 2x_i - 2p_4}{2} = \frac{x_i + x_4}{2},$$

por tanto, $S_q(m_i) = \dfrac{x_i + x_4}{2}$ para $i = 1,2,3$, obteniendo lo deseado (ver *Figura 4*).

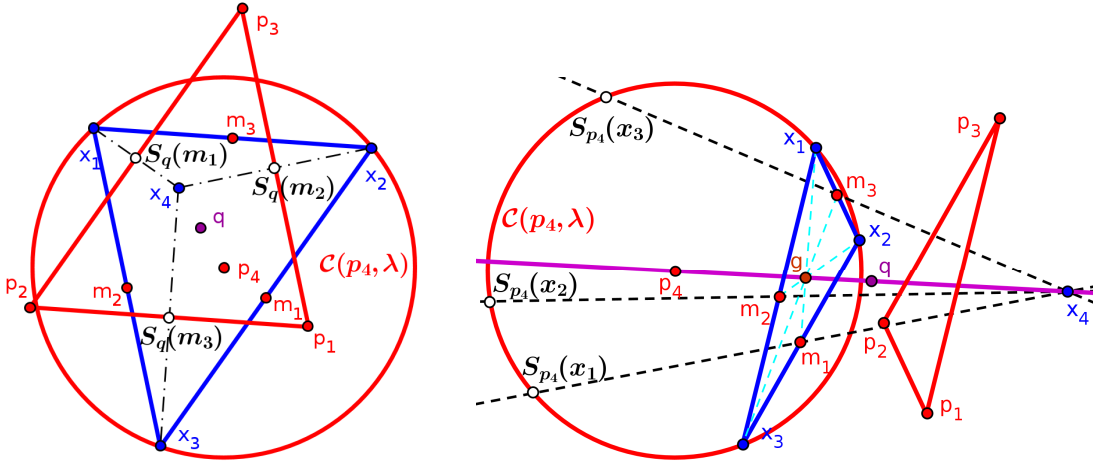

Figura 4: *Demostración del Lema 3.2*

2. Es claro que $S_q(x_i) = 2q - x_i$ para $i = 1,2,3$. Por el ítem 1 del Teorema 2.1, se tiene que

$$S_q(x_i) = x_1 + x_2 + x_3 - p_4 - x_i = x_j + x_k - p_4 = S_{m_i}(p_4).$$

para $\{i, j, k\} = \{1,2,3\}$.

3. Sea $g$ el baricentro del $\triangle x_1 x_2 x_3$. Notemos que

$$H_{g,-2}(q) = 3g - 2q = 3\left(\frac{x_1 + x_2 + x_3}{3}\right) - 2\left(\frac{x_1 + x_2 + x_3 - p_4}{2}\right) = p_4,$$

lo que dice que el segmento $[g,q]$ está contenido $2$ veces en el segmento $[p_4,g]$. Luego, por hipótesis, se tiene que $S_q(x_4)=p_4$, es decir, el segmento $[q,x_4]$ está contenido $2$ veces en el segmento $[p_4,x_4]$ y, por tanto, $q$ y $p_4$ son conjugados armónicos de los puntos $x_4$ y $g$ (ver *Figura 4*).

4. Dado que $x_4 = x_1 + x_2 + x_3 - 2p_4$, entonces

$$x_4 = 2\left(\frac{x_k+x_j}{2}\right)+(-1)(2p_4-x_i)$$

para $\{i,j,k\}=\{1,2,3\}$, y por tanto $x_4 \in \langle S_{p_4}(x_i), m_i \rangle$ para $i=1,2,3$. □

El siguiente teorema nos habla un poco sobre la relación de los puntos de Poncelet con la circunferencia de Feuerbach de un triángulo dado.

**Teorema 3.2:** Sea el $\Delta x_1 x_2 x_3$ en un plano de Minkowski. Sea $m_i$ el punto medio del segmento $[x_j, x_k]$ para $\{i,j,k\}=\{1,2,3\}$. Sean $C(p_4, \lambda)$ la circunferencia circunscrita del triángulo y $q$ el punto de simetría del triángulo dado y su $p_4$–antitriángulo, entonces:

1. Los puntos de Poncelet del $\Delta x_1 x_2 x_3$ y su $p_4$–antitriángulo están en la circunferencia de Feuerbach $C(q, \tfrac{\lambda}{2})$ del $\Delta x_1 x_2 x_3$.
2. Las circunferencias $C(m_i, \tfrac{\lambda}{2})$, para $i=1,2,3$, se intersectan en el punto $q$.

*Demostración:* Sea $x_4 \in H(\Delta x_1 x_2 x_3)$ asociado con $p_4$ y sea $w \in C(p_4, \lambda)$. Definamos $m_w = \dfrac{x_4+w}{2}$ (ver *Figura 5*) el punto medio del segmento $[x_4, w]$. Es claro que $m_w$ es un punto de Poncelet pues, por el ítem 1 del Corolario 2.1, se tiene que

$$\|q - m_w\| = \left\|\frac{x_1+x_2+x_3-p_4}{2}-\left(\frac{x_4+w}{2}\right)\right\| = \left\|\frac{p_4-w}{2}\right\| = \frac{\lambda}{2}$$

y, por tanto, los puntos de Poncelet del $\Delta x_1 x_2 x_3$ están en $C(q, \tfrac{\lambda}{2})$.

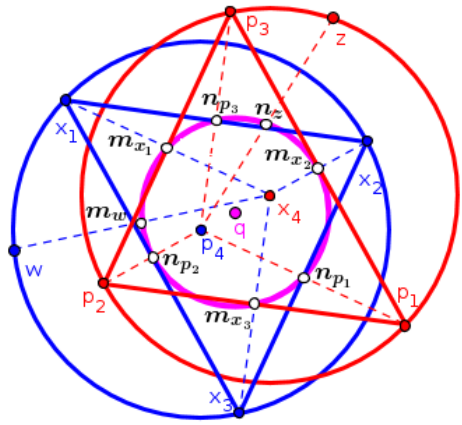

Figura 5: *Puntos de Poncelet del* $\triangle x_1 x_2 x_3$ *y su* $p_4-$*antitriángulo*.

Notemos que $x_4$ y $p_4$ son circuncentro y $C-$ortocentro del antitriángulo, respectivamente. Dado que $S_q(p_4) = x_4$, entonces $S_q(C(p_4, \lambda)) = C(x_4, \lambda)$ (ver *Figura 5*). Sea $z \in C(x_4, \lambda)$, existe $r \in C(p_4, \lambda)$ tal que $z = S_q(r)$. Denotemos por $n_z$ al punto medio del segmento $[p_4, z]$, entonces

$$\|q - n_z\| = \left\|q - \left(\frac{p_4 + z}{2}\right)\right\| = \left\|q - \left(\frac{p_4 + (2q - r)}{2}\right)\right\| = \left\|\frac{r - p_4}{2}\right\| = \frac{\lambda}{2},$$

mostrando que los puntos de Poncelet del antitriángulo están en $C(q, \frac{\lambda}{2})$.

Para mostrar que las circunferencias $C(m_i, \frac{\lambda}{2})$ se intersectan en el punto $q$, para $i = 1, 2, 3$, basta ver $\|q - m_i\| = \frac{\lambda}{2}$ $. Así,

$$\|q - m_i\| = \left\|\frac{x_1 + x_2 + x_3 - p_4}{2} - \left(\frac{x_j + x_k}{2}\right)\right\| = \left\|\frac{x_i - p_4}{2}\right\| = \frac{\lambda}{2}$$

para $i = 1, 2, 3$. □

El siguiente lema muestra una serie de sistemas $C-$ortocéntricos formados por puntos relacionados con un triángulo, tales como: puntos de Euler, circuncentro, $C-$ortocentros, puntos medios de los lados del triángulo, etc. Además, muestra algunas relaciones de paralelismo que se cumplen en el sistema $C-$ortocéntrico respectivo.

**Lema 3.3:** Sea el $\Delta x_1x_2x_3$ en un plano de Minkowski. Sea $p_4 \in C(\Delta x_1x_2x_3)$ y $x_4$ el $C$–ortocentro asociado a $p_4$. Sea $m_i$ el punto medio de los segmentos $[x_j, x_k]$ para $\{i,j,k\} = \{1,2,3\}$, entonces se cumple lo siguiente:

1. Existen $q_i \in C(m_j, \frac{\lambda}{2}) \cap C(m_k, \frac{\lambda}{2})$ tales que los puntos $q_1$, $q_2$, $q_3$ y el punto $q = \dfrac{x_4 + p_4}{2}$ (centro de la circunferencia de Feuerbach del $\Delta x_1x_2x_3$) forman un sistema $C$–ortocéntrico. Además, se cumplen las siguientes relaciones:

$$q_i - q_j = \frac{1}{2}(x_i - x_j) \quad \text{y} \quad q_i - q = \frac{1}{2}(x_i - x_4).$$

2. Si $d_i = \dfrac{x_i + x_4}{2}$, entonces $\{d_1, d_2, d_3, x_4\}$ es un sistema $C$–ortocéntrico y además, se cumplen las siguientes relaciones:

$$d_i - d_j = \frac{1}{2}(x_i - x_j) \quad \text{y} \quad d_i - x_4 = \frac{1}{2}(x_i - x_4).$$

3. El conjunto $\{m_1, m_2, m_3, p_4\}$ es un sistema $C$–ortocéntrico y además, se cumplen las siguientes relaciones:

$$m_i - m_j = \frac{1}{2}(p_i - p_j) \quad \text{y} \quad m_i - p_4 = \frac{1}{2}(p_i - p_4).$$

*Demostración:* Definamos los puntos $q_i = H_{p_4, \frac{1}{2}}(x_i)$ para $i = 1,2,3$. Así, $C(p_4, \frac{\lambda}{2})$ es una circunferencia circunscrita del $\Delta q_1q_2q_3$ y, por tanto, el punto $q' = H_{p_4, \frac{1}{2}}(q) = \dfrac{q + p_4}{2}$ es el centro de la circunferencia de Feuerbach del $\Delta q_1q_2q_3$ y $S_{q'}(p_4) = q$ es un $C$–ortocentro del mismo triangulo. Por tanto, $\{q_1, q_2, q_3, q\}$ es un sistema $C$–ortocéntrico (ver *Figura 6*).

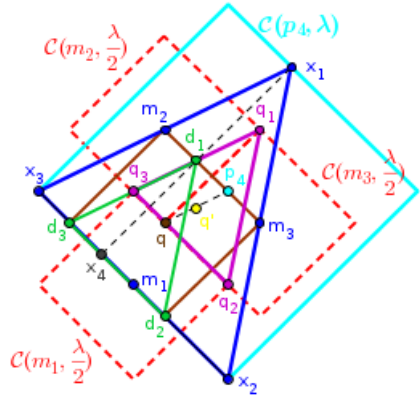

Figura 6: *Demostración del ítem 1, Lema 3.3*

Por otro lado, $q_i \in C(m_j, \frac{\lambda}{2}) \cap C(m_k, \frac{\lambda}{2})$, lo cual se puede verificar viendo que

$$\|q_i - m_j\| = \left\|\frac{p_4 + x_i}{2} - \left(\frac{x_i + x_k}{2}\right)\right\| = \frac{\|p_4 - x_k\|}{2} = \frac{\lambda}{2},$$

para $\{i, j, k\} = \{1,2,3\}$. Usando el ítem 1 del Teorema 2.1 y el ítem 3 del Corolario 3.1, se tiene que

$$q_i - q = \frac{p_4 + x_i}{2} - \left(\frac{x_1 + x_2 + x_3 - p_4}{2}\right) = p_4 - \left(\frac{x_j + x_k}{2}\right) = \frac{1}{2}(x_i - x_4),$$

para $\{i, j, k\} = \{1,2,3\}$. También,

$$q_i - q_j = H_{p_4, \frac{1}{2}}(x_i) - H_{p_4, \frac{1}{2}}(x_j) = \frac{x_i - x_j}{2},$$

y por el ítem 2 del Teorema 2.1, entonces

$$q_i - q_j = \frac{1}{2}(p_i - p_j),$$

verificando así el ítem 1 del lema.

Por otra parte, $d_i = H_{x_4, \frac{1}{2}}(x_i)$ para $i = 1,2,3$ y por tanto, $H_{x_4, \frac{1}{2}}(\Delta x_1 x_2 x_3) = \Delta d_1 d_2 d_3$, de manera que $H_{x_4, \frac{1}{2}}(x_4) = x_4$ es $C$–ortocentro del $\Delta d_1 d_2 d_3$ y así, $\{d_1, d_2, d_3, x_4\}$ es un sistema $C$–ortocéntrico, obteniendo la validez de la primera parte del ítem 2 del lema. Las relaciones faltantes se deducen usando un razonamiento similar al usado en el ítem anterior. Por último, para deducir la validez del ítem 3 se emplea el mismo razonamiento usado en el ítem 2, pero tomando en cuenta que $m_i = H_{p_4, \frac{1}{2}}(p_i)$ para $i = 1,2,3$. □

El siguiente lema muestra una propiedad de ortogonalidad isósceles que cumple todo sistema $C$ – ortocéntrico en un plano de Minkowski.

**Teorema 3.2:** Sea $\{p_1, p_2, p_3, p_4\}$ un sistema $C$ – ortocéntrico en un plano de Minkowski, entonces

$$(p_i - p_j) \perp_I (p_k - p_l)$$

para $\{i, j, k, l\} = \{1,2,3,4\}$.

**Demostración:** Sea $x_l$ un circuncentro del $\Delta p_i p_j p_k$, para $\{i, j, k, l\} = \{1,2,3,4\}$, entonces $\|x_l - p_i\| = \|x_l - p_j\|$. Por el ítem 3 del Teorema 2.1, dado que $x_l - p_i = p_j - x_k$ para $\{i, j, k, l\} = \{1,2,3,4\}$, se tiene que $p_i x_j p_j x_k$ es un paralelogramo y por tanto, $(x_l - x_k) \perp_I (p_i - p_j)$ pues

$$\|(x_l - x_k) + (p_i - p_j)\| = \|(x_l - p_j) + (p_i - x_k)\| = 2\|x_l - p_j\|$$

Y

$$\|(x_l - x_k) + (p_i - p_j)\| = \|(x_l - p_i) + (p_i - x_k)\| = 2\|x_l - p_i\|.$$

Por el ítem 2 del Teorema 2.1, se tiene que $x_l - x_k = p_k - p_l$ y por tanto, $(p_k - p_l) \perp_I (p_i - p_j)$. □

El siguiente lema muestra una relación intrínseca entre la circunferencia de Feuerbach de un triángulo y las circunferencias circunscritas del triángulo dado y su antitriángulo.

**Lema 3.4:** Sean el $\Delta x_1 x_2 x_3$ en un plano de Minkowski, $C(p_4, \lambda)$ una circunferencia circunscrita del $\Delta x_1 x_2 x_3$, $\Delta d_1 d_2 d_3$ el $p_4$ – antitriángulo del triángulo dado y $q$ el punto de simetría de estos, entonces $C(p_4, \lambda)$ es circunferencia circunscrita del $\Delta d_1 d_2 d_3$ si y solo si $q = p_4$.

*Demostración:* Por hipótesis se tiene que $d_i = S_{m_i}(p_4)$ para $i = 1,2,3$. Supongamos que los triángulos $\Delta x_1 x_2 x_3$ y $\Delta d_1 d_2 d_3$ tienen la misma circunferencia circunscrita $C(p_4, \lambda)$, entonces $x_4 = p_4$ (ver *Figura 7*). Por el ítem 3 del Teorema 3.1, obtenemos que $p_4 = q$.

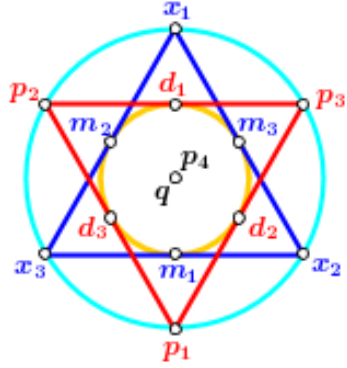

Figura 7: *Demostración del Lema 3.4*

Recíprocamente, si $p_4 = q$, se tiene que $\|p_4 - m_i\| = \dfrac{\lambda}{2}$ para $i = 1,2,3$. Por tanto,

$$\lambda = \|2p_4 - 2m_i\| = \|p_4 - (2m_i - p_4)\| = \|p_4 - S_{m_i}(p_4)\|.$$

Así, $C(p_4, \lambda)$ es una circunferencia circunscrita del $\Delta d_1 d_2 d_3$. □

El siguiente lema muestra ciertas propiedades de ortogonalidad cordal que cumple todo sistema $C$–ortocéntrico en un plano de Minkowski.

**Teorema 3.4:** Sean $M$ un plano de Minkowski, $\Delta p_1 p_2 p_3$ un triángulo de vértices distintos con circuncentro $p$, y $C$–ortocentro $h$ asociado a $p$. Sean $m_1$, $m_2$, $m_3$ los puntos medios de los segmentos $[p_2, p_3]$, $[p_3, p_1]$ y $[p_1, p_2]$, respectivamente, y sea $C$ (circunferencia unitaria) la circunferencia de Feuerbach del $\Delta p_1 p_2 p_3$, entonces se cumple:

1. Si $C \cap \langle p_i, p_j \rangle = \{m_k, v_k\}$, para $\{i, j, k\} = \{1,2,3\}$, y $C \cap [p_k, h] = \{u_k\}$, $k = 1,2,3$, entonces $[m_k, v_k] \perp_C [v_k, u_k]$.
2. Si $u_k = \dfrac{p_k + h}{2}$ para $k = 1,2,3$, entonces $[m_i, -m_j] \perp_C [h, u_k]$ para $\{i, j, k\} = \{1,2,3\}$, y $[m_i, m_j] \perp_C [u_i, -u_j]$ para $\{i, j\} \subset \{1,2,3\}$.

*Demostración:* 1. Puesto que $h$ es $C$–ortocentro del $\Delta\, p_1 p_2 p_3$, se tiene un sistema $C$–ortocéntrico $\{p_1, p_2, p_3, h\}$. Como $p$ es circuncentro de dicho triángulo, entonces por el Teorema 3.1, se tiene que el centro de la circunferencia de Feuerbach del triángulo está dado por $\dfrac{p+h}{2}$, y como por hipótesis esta circunferencia coincide con la unitaria, se tiene que $O = \dfrac{p+h}{2}$, lo que implica que $p = -h$.

Por el ítem 3 del Corolario 2.1, se tiene que $h = p_1 + p_2 + p_3 - 2p$ y así, $O = p_1 + p_2 + p_3 + h$ de manera que:

$$-\frac{p_1 + p_2}{2} = \frac{p_3 + h}{2}, \qquad -\frac{p_1 + p_3}{2} = \frac{p_2 + h}{2}, \qquad -\frac{p_2 + p_3}{2} = \frac{p_1 + h}{2}.$$

Como $u_k$ es el punto medio del segmento $[p_k, h]$, se tiene que $u_k = -m_k$ para $k = 1,2,3$ y, por tanto, son puntos opuestos de $C$. Dado que $m_k$ y $v_k$ están en un mismo lado del triángulo (ver *Figura 8*), entonces $v_k \neq -m_k$. De manera similar $v_k \neq -u_k$, pues de ser iguales se tendría que $v_k = m_k$, lo que contradice la hipótesis. Luego, por la Proposición 2.1, se tiene que $[m_k, v_k] \perp_C [v_k, u_k]$ para $k = 1,2,3$.

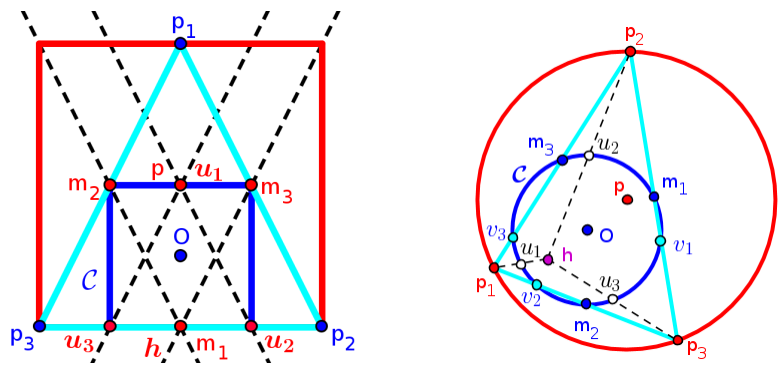

Figura 8: *Demostración del Teorema 3.4*

2. Por el ítem 1 del Teorema 2.1, se tiene que $O = \dfrac{p+h}{2}$. Por el ítem 3 del Corolario 2.1, $h = -p_1 - p_2 - p_3$. Como $u_k = \dfrac{p_k + h}{2}$ entonces

$$-h - u_k = \frac{p_k - h}{2} = \frac{(p_k + p_i) + (p_j + p_k)}{2} = m_i + m_j,$$

para $\{i,j,k\} = \{1,2,3\}$. Así, el vector director de la recta $\langle m_i, -m_j \rangle$ es múltiplo escalar del vector director de la recta $\langle -h, u_k \rangle$, es decir, $[m_i, -m_j]$ es paralelo a $[-h, u_k]$ para $\{i,j,k\} = \{1,2,3\}$ y por tanto,
$[m_i, -m_j] \perp_C [h, u_k]$.

Notemos que

$$u_i - u_j = \frac{h + p_i}{2} - \left(\frac{h + p_j}{2}\right) = \frac{p_i - p_j}{2} = \frac{(p_i + p_k) - (p_j + p_k)}{2} = m_i - m_j.$$

De manera que la recta $\langle m_i, m_j \rangle$ es paralela a la recta $\langle u_i, u_j \rangle$ para $\{i,j\} \subset \{1,2,3\}$ y por tanto, $[m_i, m_j] \perp_C [u_i, -u_j]$ para $\{i,j\} \subset \{1,2,3\}$. □

El siguiente teorema dice que la imagen homotética de un sistema $C$ – ortocéntrico, en un plano de Minkowski, es también un sistema $C$ – ortocéntrico.

**Teorema 3.5:** La imagen homotética de un sistema $C$ – ortocéntrico es un sistema $C$ – ortocéntrico.

*Demostración:* Sea $\{p_1, p_2, p_3, p_4\}$ un sistema $C$ – ortocéntrico, entonces existe $x_4 \in C(\Delta p_1 p_2 p_3)$ tal que $p_4 = p_1 + p_2 + p_3 - 2x_4$. Sean $H_{w,k}$ una homotecia de centro en un punto $w$ y razón $k$, $h_i = H_{w,k}(p_i)$ para $i = 1, 2, 3$, y $y_4 = H_{w,k}(x_4)$. Claramente $y_4 \in C(\Delta h_1 h_2 h_3)$ y

$$h_1 + h_2 + h_3 - 2y_4 = ((1-k)w + kp_1) + ((1-k)w + kp_2) + ((1-k)w + p_3) - 2((1-k)w + kx_4) =$$
$$= (1-k)w + p_1 + p_2 + p_3 - 2x_4 = (1-k)w + p_4 = h_4$$

lo que completa la demostración. □

El siguiente corolario dice que dado un triángulo, con $C$ – ortocentro asociado a un circuncentro del mismo, el conjunto formado por los baricentros de los cuatro triángulos que se obtienen con los vértices del triángulo y el $C$ – ortocentro dado, es un sistema $C$ – ortocéntrico.

**Corolario 3.3:** Sea $\{p_1, p_2, p_3, p_4\}$ un sistema $C$ – ortocéntrico en un plano de Minkowski arbitrario. Sea $g_i$ el baricentro y $x_i$ circuncentro del $\Delta p_j p_k p_l$ para $\{i,j,k,l\} = \{1,2,3,4\}$. Entonces $\{g_1, g_2, g_3, g_4\}$ es un sistema $C$ – ortocéntrico y

$$g_i - g_j = \frac{1}{3}(p_j - p_i) = \frac{1}{3}(x_i - x_j).$$

*Demostración:* Dado que $g_i = \dfrac{p_j + p_k + p_l}{3}$ para $\{i, j, k, l\} = \{1,2,3,4\}$, entonces

$$g_i - g_j = \frac{p_j + p_k + p_l}{3} - \left(\frac{p_i + p_k + p_l}{3}\right) = \frac{p_j - p_i}{3} = \frac{x_i - x_j}{3}.$$

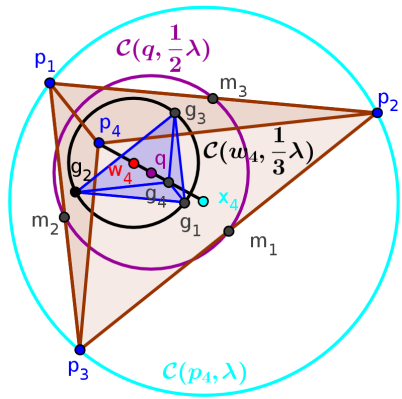

Figura 9: *Sistema $C$ – ortocéntrico de baricentros.*

Sea $\{p_1, p_2, p_3, p_4\}$ un sistema $C$ – ortocéntrico, entonces existe $x_4 \in C(\Delta p_1 p_2 p_3)$. Sea $q$ el punto de simetría del $\Delta p_1 p_2 p_3$ y su $x_4$ –-antitriángulo. Veamos que $g_i = H_{q,-\frac{1}{3}}(p_i)$ para $i = 1,2,3,4$ (ver *Figura 9*):

$$H_{q,-\frac{1}{3}}(p_i) = \frac{4}{3}q - \frac{1}{3}p_i = \frac{4}{3}\left(\frac{x_i + p_i}{2}\right) - \frac{1}{3}p_i = \frac{1}{3}p_i + \frac{2}{3}x_i =$$
$$= \frac{1}{3}(p_j + p_k + p_l - 2x_i) + \frac{2}{3}x_i = g_i$$

pues por el ítem 3 del Corolario 3.1, $p_i = p_j + p_k + p_l - 2x_i$ para $\{i, j, k, l\} = \{1,2,3,4\}$. Luego, por el Teorema 3.5 se tiene que $\{g_1, g_2, g_3, g_4\}$ es un sistema $C$ – ortocéntrico. □

**Referencias.**